\documentstyle[a4,leqno,amsfonts,12pt]{article}
\newcommand{\C}{{\bf C}}

\newcommand{\OC}{\overline{{\bf C}}}
\newcommand{\R}{{\bf R}}

\newcommand{\N}{{\bf N}}
\newcommand{\M}{{\bf M}}

\newcommand{\D}{{\bf D}}
\newcommand{\de}{\delta}

\newcommand{\mb}{\mbox}
\newcommand{\di}{\mb{dist}}

\newcommand{\beq}{\begin{equation}}
\newcommand{\eeq}{\end{equation}}

\newcommand{\ov}{\overline}

\newcommand{\Om}{\Omega}

\newcommand{\z}{\zeta}

\newcommand{\kap}{\mb{cap}}

\newtheorem{lem}{Lemma}

\newcommand{\ueberschrift}{\bigskip\goodbreak\noindent\bigskip}
\newcounter{theabsatz}
\newcommand{\absatz}[1]{\stepcounter{theabsatz} \ueberschrift
                           {\large \bf \arabic{theabsatz}. {#1}} \setcounter{equation}{0}}

\begin{document}
\mathsurround=2pt

\begin{center}
{\large \bf ON THE TOTIK-WIDOM PROPERTY FOR A QUASIDISK}\\[2ex]
V. Andrievskii, F. Nazarov\\[2ex]

\end{center}

\begin{abstract}

Let $K$ be a quasidisk on the complex plane. We construct a sequence
of monic polynomials $p_n=p_n(\cdot,K)$ with zeros on $K$ such that
\newline $||p_n||_K\le O(1)\kap(K)^n$ as $n\to\infty.$

\end{abstract}




\absatz{Introduction and Main Result}

Let $K\subset\C$ be a {\it continuum}, i.e., a compact set in the complex plane $\C$
with a simply connected complement $\Om:=\ov{\C}\setminus K$, where
$\ov{\C}:=\C\cup\{\infty\}.$ We assume that $\kap(K)>0$,
where$\kap(S)$ denotes the logarithmic capacity of a compact set
$S\subset\C$ (see \cite{ran, saftot}).
Let $\M_n$ be the class of all monic polynomials of degree $n
\in \N:=\{ 1,2,\ldots\}$.
Denote by $\tau_n(z)=\tau_n(z,K),n\in\N$
the $n$-th {\it Chebyshev polynomial} (associated with $K$), i.e.,
$\tau_n\in\M_n$ is the (unique) polynomial which minimizes the uniform norm
$||\tau_n||_K:=\sup_{z\in K}|\tau_n(z)|$ among all monic polynomials of the same
degree.
Let $\tilde{\tau}_n(z)=\tilde{\tau}_n(z, K)$ be
the $n$-th  Chebyshev polynomial with zeros on $K$ (see \cite[p. 155]{ran}
and \cite[Chapter 1, \S 3]{smileb}). 

It is well-known (see, for example, \cite[Theorem 5.5.4 and
Corollary 5.5.5]{ran} and \cite[Chapter 1, \S 3]{smileb})  that
$$
||\tilde{\tau}_n||_K\ge||\tau_n||_K\ge\kap(K)^n\quad\mb{and}\quad
\lim_{n\to\infty}||\tilde{\tau}_n||_K^{1/n}=\lim_{n\to\infty}||\tau_n||_K^{1/n}=\kap(K).
$$
Let
$$
w_n(K):=\frac{||\tau_n||_K}{\kap(K)^n}\quad\mb{and}\quad \tilde{w}_n(K):=\frac{||\tilde{\tau}_n||_K}{\kap(K)^n}
$$
be the {\it Widom factors}. The estimate of $w_n(K)$ from above has
attracted considerable attention.
For the complete survey of results concerning  this problem and further citations, see
\cite{csz}, \cite{smileb}-\cite{wid},  \cite{and17}. At the same time, $\tilde{\tau}_n$ and their zeros play
significant role in interpolation theory (cf. \cite[Chapter 1, \S 3]{smileb}).

Recently, Simon \cite{sim} asked the question whether the closed domain $K$
bounded by the Koch snowflake  obeys 
a {\it Totik-Widom bound}, i.e.,
\beq
\label{1.tw}
w_n(K)=O(1)\quad \mb{as }n\to\infty?
\eeq
To answer this question,
we will consider the case where $K=\ov{G}$ is the closure of a Jordan domain $G$ bounded
by a quasiconformal curve $L:=\partial G$ (a {\it quasidisk} for short),
see \cite{ahl} or \cite{ger}.
Note that by the Ahlfors criterion (see \cite[p. 100]{lehvir}) the Koch snowflake
is quasiconformal.

{\bf Theorem.} {\it
Let $K$
be a quasidisk.
 Then}
\beq\label{1.5}
\widetilde{w}_n(K)=O(1)\quad \mb{as }n\to\infty.
\eeq

In particular, any quasidisk
obeys the Totik-Widom bound (\ref{1.tw}), 
see also \cite[Theorem 2]{and17} where quite 
different method for the estimate of $w_n(K)$ is used.

\absatz{Auxiliary Constructions}

Till the end of this section $K$ is a continuum with $\kap(K)>0$. 
Below we use the notation
$$
 \D:=\{ z:|z|<1\},\quad \D^*:=\ov{\C}\setminus\ov{\D},
$$
$$
\di (S_1,S_2):=\inf_{z_1\in S_1,z_2\in S_2}|z_2-z_1|,\quad S_1,S_2\subset\C.
$$
Let $\Phi:\Om\to\D^*$  be the
Riemann conformal mapping normalized by the conditions
$$
\Phi(\infty)=\infty,\quad
\Phi'(\infty):=\lim_{z\to\infty}\frac{\Phi(z}{z}>0,
$$
$\Psi:=\Phi^{-1}$,
 and let
$$
K_s:=\{z\in \Om:|\Phi(z)|=1+s\},\quad s>0.
$$
For $0<s<1$,  denote by $\mu_s$ the equilibrium measure for
$K_s$ (see \cite{ran}, \cite{saftot}). 
Let $I'=I'(s,\theta_0,\theta_1):=\{w=(1+s)e^{i\theta}:\theta_0\le\theta\le\theta_1\}, 
0<\theta_1-\theta_0<s$
and $I:=\Psi(I'),\z_0:=\Psi((1+s)e^{i\theta_0}).$
For an arbitrary (but fixed) $q\in\N\setminus\{1\}$ and
$$
m_l:=\frac{1}{\mu_s(I)}\int_I(\xi-\z_0)^ld\mu_s(\xi),\quad l=1,\ldots,q,
$$
consider the system of equations (with respect to $r_k,k=1,\ldots,q$)
$$
\sum_{k=1}^qr_k^l=qm_l=:\tilde{m}_l, \quad l=1,\ldots,q.
$$
We interpret $r_k=r_k(I,q)$ as the roots of the polynomial
$
z^q+a_{q-1}z^{q-1}+\ldots+a_0
$
whose coefficients satisfy Newton's identities (see, for example \cite{kal})
\beq\label{2.2}
\tilde{m}_l+a_{q-1}\tilde{m}_{l-1}+\ldots+a_{q-l+1}\tilde{m}_1=-la_{q-l},
\quad l=1,\ldots, q.
\eeq
Since
$|m_l|\le d^l$, where $d:=$diam$I$, according to (\ref{2.2}) 
 we obtain
\beq\label{2.3}
|a_{q-l}|\le q^ld^l,\quad l=1,\ldots,q
\eeq
(the proof is by induction on $l$).

Furthermore, if $r_k\neq0$ then by (\ref{2.3}),
\begin{eqnarray}
1&\le&\frac{|a_{q-1}|}{|r_k|}+
\frac{|a_{q-2}|}{|r_k|^2}+\ldots+ \frac{|a_{0}|}{|r_k|^q}\nonumber\\
\label{2p.1}
&\le& \frac{qd}{|r_k|}+
\frac{q^2d^2}{|r_k|^2}+\ldots+ \frac{q^qd^q}{|r_k|^q} ,
\end{eqnarray}
from which we derive
\beq\label{2.4}
|r_k|\le 2q d,\quad k=1,\ldots,q.
\eeq
Indeed, supposing contrary to our claim that $|r_k|>2qd$ and
using (\ref{2p.1})
we get the contradiction
$$
1\le \frac{qd}{|r_k|}
\sum_{j=0}^\infty\left(\frac{qd}{|r_k|}\right)^j=
\frac{qd}{|r_k|}\left(1-\frac{qd}{|r_k|}\right)^{-1}<\frac{2qd}{|r_k|}<1.
$$
For  $m\in\N\setminus\{1\}$, $0<s<1$, and
 $j=1,\ldots, m$, consider
$$
I_j':=\left\{w=(1+s)e^{i\theta}:\frac{2\pi(j-1)}{m}\le\theta\le\frac{2\pi j}{m}\right\},
$$
$$
\xi_j:=\Psi\left( (1+s)e^{2\pi(j-1)/m}\right),\quad
\xi_{m+1}:=\xi_1,\quad I_j:=\Psi(I_j').
$$
Let $\lfloor x \rfloor$ denote the integer part of $x\in\R$ and let $|S|$
be the linear measure (length) of a (Borel) set $S\subset\C$.
\begin{lem}\label{lem2.1}
For $s=c q/m$, where $c=480\pi ,m,q\in\N,m>m_0:=\lfloor cq\rfloor +1$
 and $j=1,\ldots,m$, we have
\beq\label{2.5}
\frac{{\em \di}(I_j,K)}{2000\pi q}\le |\xi_{j+1}-\xi_j|\le \mb{\em diam }I_j\le
|I_j|\le \frac{{\em \di}(I_j,K)}{10q} .
\eeq
\end{lem}
{\bf Proof.}
By \cite[p. 23, Lemma 2.3]{andbla}, which is an immediate consequence of
Koebe's one-quarter theorem, for $\tau \in\D^*\setminus\{\infty\}$
and $\z:=\Psi(\tau)$ we have
\beq\label{n2.1}
\frac{1}{4}\frac{\di(\z,K)}{|\tau|-1}\le|\Psi'(\tau)|\le
4\frac{\di(\z,K)}{|\tau|-1}.
\eeq
Moreover, if $|\eta-\tau|\le(|\tau|-1)/2$ and $\xi:=\Psi(\eta)$, then
\beq\label{n2.2}
\frac{1}{16}\frac{|\eta-\tau|}{|\tau|-1}\le\frac{|\xi-\z|}{\di(\z,K)}\le
16\frac{|\eta-\tau|}{|\tau|-1}.
\eeq
Let $\tau_j:=\Phi(\xi_j).$ Since for $\tau\in I'_j$
$$
|\tau-\tau_j|\le|\tau_{j+1}-\tau_j|\le(1+s)\frac{2\pi}{m}\le\frac{4\pi}{m}=\frac{4\pi s}{cq},
$$
 by (\ref{n2.2})  we obtain for $\xi=\Psi(\tau)\in I_j$
$$
\frac{|\xi-\xi_j|}{\di(\xi_j,K)}\le16\frac{4\pi s}{cqs}=\frac{64\pi}{cq}<\frac{1}{2}.
$$
Therefore, by virtue of (\ref{n2.1}), we have
\begin{eqnarray*}
|I_j|&=&\int_{I'_j}|\Psi'(\tau)||d\tau|\\
&\le&\frac{4(1+s)}{s}\int_{2\pi(j-1)/m}^{2\pi j/m}\di(\Psi((1+s)e^{i\theta}),K)d\theta\\
&\le& \frac{24\pi}{sm}\di(\xi_j,K)\le\frac{48\pi}{cq}\di(I_j,K)=\frac{\di(I_j,K)}{10q},
\end{eqnarray*}
which implies the last inequality in (\ref{2.5}).

The first inequality in (\ref{2.5}) follows immediately from (\ref{n2.2})
as follows
\begin{eqnarray*}
\frac{|\xi_{j+1}-\xi_j|}{\di(I_j,K)}&\ge& \frac{|\xi_{j+1}-\xi_j|}{\di(\xi_j,K)}
\ge \frac{1}{16}\frac{2}{\pi}
\frac{2\pi(1+s)}{m}\frac{1}{s}\\
&\ge& \frac{1}{4cq}>\frac{1}{2000\pi q}.
\end{eqnarray*}

 \hfill$\Box$

\begin{lem}\label{lem2.2}
For any $m, q\in\N, m>cq$, where $c=480\pi$,
and $n=qm$ there exists a polynomial $p_n\in \M_n$ with all zeros in the interior
of $K_{5cq/m}$ such that
\beq\label{1.1}
||p_n||_K\le c_1\exp\left( c_2\left|\left|\int_{K_{cq/m}}\frac{d(\z,K)^{q}|d\z|}{|
\z-\cdot|^{q+1}}\right|\right|_{\partial K}\right)\mb{\em {cap}} (K)^n,
\eeq
where $c_j=c_j(q)>0$, j=1,2.
\end{lem}
 {\bf Proof.}
We construct the  points $\z_j^k,  j =1,\ldots, m ,  k=1,\ldots, q$ as follows. 
For $l=1,\ldots, q$, let
$$
 m_{j,l}:=\frac{1}{\mu_s(I_j)}
\int_{I_j}(\xi-\xi_j)^ld\mu_s(\xi),\quad s:=\frac{cq}{m}.
$$
The system
$$
\sum_{k=1}^qr_{j,k}^l=qm_{j,l}
$$
has solutions $r_{j,k}$ satisfying, according to (\ref{2.4}), the inequality
\beq\label{2n.1}
|r_{j,k}|\le 2qd_j,\quad d_j:=\mb{diam }I_j.
\eeq
Let $\z_j^k:=\xi_j+r_{j,k}.$
Since by virtue of (\ref{2.5}) and (\ref{2n.1})
$$
\frac{|\z_j^k-\xi_j|}{\di(\xi_j,K)}\le \frac{|r_{j,k}|}{\di(I_j,K)}\le \frac{1}{5},
$$
(\ref{n2.2}) implies
$$
|\Phi(\z_j^k)-\Phi(\xi_j)|\le\frac{16s}{d(\xi_j,K)}|\z_j^k-\xi_j|\le\frac{16s}{5}
$$
which further yields
$
|\Phi(\z_j^k)|<1+5s,
$
i.e., all $\z_j^k$ belong to the interior of $K_{5s}$.

Consider the polynomial
$$
p_n(z):=\prod_{j=1}^m\prod_{k=1}^{q}(z-\z_j^k),\quad n=qm.
$$
Since $\mu_s(I_j)=1/m$ and for $z\in \partial K$,
$$
\log\kap(K_s)=\int_{K_s}\log|z-\xi|d\mu_s(\xi)
$$
(see \cite[pp. 59, 127]{ran}), we have
\begin{eqnarray}
&&\log|p_n(z)|-n\log\kap(K_s)\nonumber\\
&=&
\sum_{j=1}^m\sum_{k=1}^q\left(\log|z-\z_j^k|-m
\int_{I_j}\log|z-\xi|d\mu_s(\xi)\right)\nonumber\\
&=&
m\sum_{j=1}^m\sum_{k=1}^q\int_{I_j}\log\left|\frac{z-\z_j^k}{z-\xi}\right|
d\mu_s(\xi)\nonumber\\
\label{2.pu1}
&\le&
m\sum_{j=1}^m\left|\sum_{k=1}^q\int_{I_j}\log\frac{z-\z_j^k}{z-\xi}
d\mu_s(\xi)\right|.
\end{eqnarray}
Next,  by virtue of Lemma \ref{lem2.1}, (\ref{2n.1}), 
and Taylor's Theorem \cite[pp. 125-126]{ahl1}
for $\xi\in I_j$ and $z\in \partial K$,
\begin{eqnarray}
\log\frac{z-\z_j^k}{z-\xi}&=&
\log\left(1-\frac{\z_j^k-\xi_j}{z-\xi_j}\right)-
\log\left( 1-\frac{\xi-\xi_j}{z-\xi_j}\right) \nonumber\\
&=&
\label{2v.1}
\sum_{l=1}^q\frac{1}{l}\left( \left( \frac{\xi-\xi_j}{z-\xi_j}\right)^l
-\left(\frac{\z_j^k-\xi_j}{z-\xi_j}\right)^l\right)
+B_j(z),
\end{eqnarray}
where
$$
|B_j(z)|=|B_j(z,\xi,\{\z_j^k\}_{k=1}^q)|\le c_3
\left(\frac{d_j}{\di(z,I_j)}\right)^{q+1},\quad c_3=c_3(q).
$$
Since
$
\kap(K_s)=(1+s)\kap(K)
$
and
\begin{eqnarray*}
\sum_{k=1}^q\int_{I_j}\left( (\xi-\xi_j)^l-(\xi_j^k-\xi_j)^l\right) d\mu_s(\xi)
&=& \sum_{k=1}^q\left(\frac{1}{m} m_{j,l}-\frac{1}{m} r^l_{j,k}\right)\\
&=& \frac{1}{m}\left( qm_{j,l}-\sum_{k=1}^q r_{j,l}^l\right)=0,
\end{eqnarray*}
applying (\ref{2.pu1}) and (\ref{2v.1}) we obtain
\begin{eqnarray*}
|\log|p_n(z)|-n\log\kap(K)|&\le& n\log(1+s) +|\log|p_n(z)|-n\log\kap(K_s)|
\\
&\le& cq^2 +c_3q
\sum_{j=1}^m\left(\frac{d_j}{\di(z,I_j)}\right)^{q+1}.
\end{eqnarray*}
Furthermore, since by Lemma \ref{lem2.1}
\begin{eqnarray*}
&&\int_{K_s}\frac{\di(\z,K)^q}{|\z-z|^{q+1}}|d\z|=
\sum_{j=1}^m\int_{I_j}\frac{\di(\z,K)^q}{|\z-z|^{q+1}}|d\z|\\
&\ge&\sum_{j=1}^m\frac{(10qd_j)^q|I_j|}{(\di(z,I_j)+d_j)^{q+1}}
\ge c_4
\sum_{j=1}^m\left(\frac{d_j}{\di(z,I_j)}\right)^{q+1},
\end{eqnarray*}
where $c_4=c_4(q)$,
we obtain (\ref{1.1}).

 \hfill$\Box$

\absatz{Proof of Theorem}

Let $Q>1$ be the coefficient of quasiconformality of $L$, i.e.,
$L$ is the image of the unit circle under some $Q$-quasiconformal
mapping $F:\OC\to\OC$.
Denote by $c_5, c_6,\ldots$ positive constants
that either are absolute, or depend 
on $Q$ only.
According to
\cite[Chapter IV]{ahl} the Riemann mapping function $\Phi$
can be extended to a $Q^2$-quasiconformal homeomorphism $\Phi:\OC\to\OC$.

For $\z\in\Om$ denote by $\z_K:=\Psi(\Phi(\z)/|\Phi(\z)|)$
the ``projection" of $\z$ to $K$.
As an immediate application of \cite[p. 29, Theorem 2.7]{andbla}
(cf. \cite[p. 21, Theorem 2.6 and p. 23, Corollary 2.7]{ger})
for $\z\in\Om$ and $z\in L$ we have
$$
\frac{\di(\z,K)}{|\z-z|}\le\left|\frac{\z-\z_K}{\z-z}\right|
\le c_5\left(\frac{|\Phi(\z)|-1}{|\Phi(\z)-\Phi(z)|}\right)^{1/Q^2}.
$$
Hence, according to (\ref{n2.1}) 
for $z\in L, w=\Phi(z)$, and $0<s<1$ we obtain
\begin{eqnarray}
\int_{K_s}\frac{\di(\z,K)^q|d\z|}{|\z-z|^{q+1}}
&\le&
\frac{4}{s}\int_{|\tau|=1+s}\left(\frac{\di(\Psi(\tau),K)}{|\Psi(\tau)-
\Psi(w)|}\right)^{q+1}|d\tau|\nonumber\\\label{3.1}
&\le&\frac{4c_5^{q+1}}{s}\int_{|\tau|=1+s}\frac{s^{(q+1)/Q^2}|d\tau|}{|\tau-w|^{(q+1)/Q^2}}\le
c_6,
\end{eqnarray}
if $q=\left\lfloor2Q^2\right\rfloor$.

Therefore, if we let in (\ref{1.1}) $q:=\left\lfloor 2Q^2\right\rfloor$ then for any $n$ of the form
$n=mq$ with sufficiently large $m>m_0(Q)$ we 
have the polynomial $p_n(\cdot, K)\in\M_n$ with zeros in the interior of
 $K_{c_7/n}$ such that
$$
||p_n(\cdot, K)||_K\le c_8\kap(K)^n.
$$
Next, we introduce  auxiliary families of quasiconformal curves and mappings
as follows. 
Note that each curve
$$
L^*_\de:=\{z:|\Phi(z)|=1-\de\},\quad 0<\de<1,
$$
is $Q^2$-quasiconformal. Denote by $\Om^*_\de$ the unbounded connected component
of $\OC\setminus L^*_\de$. The Riemann conformal mapping $\Phi_\de:\Om^*_\de\to
\D^*$ with the normalization
$$
\Phi_\de(\infty)=\infty,\quad \Phi_\de'(\infty)>0
$$
can be extended to a $Q^4$-quasiconformal homeomorphism 
$\Phi_\de:\OC\to\OC$.

By \cite[p. 29, Theorem 2.7 and p. 376, Lemma 2.2]{andbla} for
$\z\in L$ and $0<\de<1/2$,
$$
\frac{\de}{c_9}\le |\Phi_\de(\z)|-1\le c_9\de.
$$
Let $c_{10}:=c_7c_9, K_n:=\C\setminus\Om^*_{c_{10}/n}, n=mq>2c_{10}$.
 Consider polynomial
$p_n(\cdot,K_n)\in\M_n$ constructed as above only using $K_n$ instead of
$K$. Since for $\z\in L$,
$$
|\Phi_{c_{10}/n}(\z)|-1\ge\frac{c_{10}}{c_9n}=\frac{c_7}{n},
$$
all zeros of $p_n(\cdot,K_n)$ belong to $K$. Moreover,
by the
Bernstein-Walsh Lemma (see \cite[p. 153]{saftot} or \cite[p. 156]{ran}) 
as well as (\ref{1.1}) and (\ref{3.1}) we have
$$
||p_n(\cdot,K_n)||_K\le \left(1+\frac{c_9c_{10}}{n}\right)^n
||p_n(\cdot,K_n)||_{K_n}\le c_{11}\kap(K)^n,
$$ 
which implies (\ref{1.5}) for $n=mq$ with $m>c_{12}$.

For an arbitrary $n\in\N$ (\ref{1.5}) follows from 
this particular case and an obvious inequality
 $$
 \tilde{w}_{n+m}(K)\le \tilde{w}_n(K)\tilde{w}_m(K),\quad n,m\in\N.
$$

 \hfill$\Box$


\begin{thebibliography}{99}



\bibitem{ahl}
L. V. Ahlfors (1966): Lectures on Quasiconformal Mappings.
Princeton, N.J.: Van Nostrand.

\bibitem{ahl1}
L. V. Ahlfors (1979): Complex Analysis. McGraw-Hill, Inc.
	
	\bibitem{and17}
 V. V. Andrievskii (2017): {\it On Chebyshev polynomials in the complex plane.}
 Acta Math. Hungar., {\bf 152(2)}:505--524.

 \bibitem{andbla}
  V. V. Andrievskii, H.- P. Blatt (2002): Discrepancy of Signed Measures and
      Polynomial Approximation.
      Berlin/New York: Springer-Verlag.


\bibitem{csz}
J. S. Christiansen, B. Simon, and M. Zinchenko (2017): {\it Asymptotics of Chebyshev polynomials,
I. Subsets of \R}, 
Invent. Math., {\bf 208}:217--245. 

\bibitem{ger}
F. W. Gehring (1982): Characteristic Properties of Quasidisks. 
Montr\'eal: Presses de l'Universit\'e de Montr\'eal.

\bibitem{kal}
D. Kalman (2000): {\it A matrix proof of Newton's identities},
Mathematics Magazine, {\bf 73(4)}:313--315.

\bibitem{lehvir}
O. Lehto, K. I. Virtanen (1973):  Quasiconformal Mappings in the Plane, 2nd
ed.,   New York: Springer-Verlag.


\bibitem{ran}
T. Ransford (1995): Potential Theory in the Complex Plane,
Cambridge: Cambridge University Press.


\bibitem{saftot}
 E. B. Saff ,  V. Totik (1997):    Logarithmic Potentials with External
Fields,  New York/Berlin: Springer-Verlag.

\bibitem{sim}
B. Simon (2017): {\it Szeg\"{o}-Widom asymptotics for Chebyshev polynomials
on subsets of $\R$}, Computational Methods and Function Theory, July
10-15, Lublin, Poland.

\bibitem{smileb}
V. I. Smirnov, N. A. Lebedev (1968):   Functions of a Complex
Variable. Constructive Theory,  Cambridge: Mass. Institute of
Technology.

\bibitem{sodyud}
M. L. Sodin, P. M. Yuditskii (1993): {\it Functions least deviating
from zero on closed subsets of the real line.} St. Petersburg Math. J.,
{\bf 4}:201--249.


\bibitem{tot12}
V. Totik (2012): {\it Chebyshev polynomials on a system of curves.}
Journal D'Analyse Math\'ematique, {\bf 118}:317--338.

\bibitem{tot13}
V. Totik (2014): {\it Chebyshev polynomials on compact sets.}
Potential Anal., {\bf 40(4)}:511–-524.

\bibitem{totvar}
V. Totik, T. Varga (2015): {\it Chebyshev and fast decreasing polynomials.}
Proc. London Math. Soc., {\bf 110(5)}:1057--1098.

\bibitem{wid}
H. Widom (1969): {\it Extremal polynomials assosiated with a system
of curves in the complex plane.} Adv. Math., {\bf 3}:127--232.

\end{thebibliography}
\end{document}